\documentclass{article}
\usepackage{amsmath}
\usepackage{extarrows} 
\usepackage{amssymb}
\usepackage[normalem]{ulem}
\usepackage{arxiv}

\usepackage[utf8]{inputenc} 
\usepackage[T1]{fontenc}    
\usepackage{hyperref}       
\usepackage{url}            
\usepackage{booktabs}       
\usepackage{amsfonts}       
\usepackage{nicefrac}       
\usepackage{microtype}      
\usepackage{amsthm}

\usepackage{hyperref} 
\usepackage{theoremref}

\title{Lyapunov stability of compact sets in locally compact metric spaces}
\author{
  Reza Hadadi\\
  Department of Electrical and Computer Engineering\\
  University of Maryland\\
  College Park, MD 20784 \\
  \texttt{rezahd@umd.edu} \\
}

\usepackage{hyperref} 
\usepackage{theoremref}

\newtheorem{theorem}{Theorem}
\newtheorem*{theorem*}{Theorem}
\newtheorem{lemma}[theorem]{Lemma}
\newtheorem{proposition}[theorem]{Proposition}
\theoremstyle{definition}
\newtheorem{definition}{Definition}

\begin{document}
\maketitle
\thispagestyle{empty}
\pagenumbering{arabic}
\begin{abstract}
    This paper provides a systematic exposition of Lyapunov stability for compact sets in locally compact metric spaces. We explore foundational concepts, including neighborhoods of compact sets, invariant sets, and the properties of dynamical systems, and establish key results on the relationships between attraction, invariance, and stability. The work explores Lyapunov stability within the context of dynamical systems, highlighting equivalent formulations and related criteria. Central to the exposition is a proof of the fundamental theorem linking Lyapunov functions to the asymptotic stability of compact sets. This expository piece consolidates results from several classical texts to provide a unified and accessible framework for understanding stability in metric spaces.
\end{abstract}
Throughout this paper assume $(X,d)$ to be a locally compact metric space and $M$ to be a non-empty compact subset of $X$, unless stated otherwise. Our standard notations that we adopt are organized in the following table:

\begin{center}
\begin{tabular}{ l| p{10cm} }
\hline
$(X,d)$ & the generic metric space,\\
$2^X$ & power set of $X$,\\
 $\mathbb{R}$ & set of real numbers,  \\ 
 $\mathbb{R}_+$ & set of non-negative real numbers,\\
 $\mathbb{R}_{-}$ & set of non-positive real numbers,\\
 $\mathbb{N}$& set of positive integers,\\
 $A^{\mathrm{o}}$& interior of the set $A$,\\
 $\partial A$& boundary of the set $A$,\\
 $\overline{A}$& closure of the set $A$,\\
 $B(x,r)$ & the open ball of radius $r>0$ centered at $x \in X$, i.e., the set $\{y \in X: d(y,x)<r\}$,\\
 $B(A,r)$ & the set $\{x \in X: d(x,A)<r\}$, where $A \subset X$ and $r>0$,\\
 $S[x,r]$ & the closed ball of radius $r\geq 0$ centered at $x \in X$, i.e., the set $\{y \in X: d(y,x)\leq r\}$,\\
 $S[A,r]$ & the set $\{x \in X: d(x,A)\leq r\}$, where $A \subset X$ and $r\geq 0$,\\
 $H(x,r)$ & the spherical hypersurface of radius $r\geq 0$ centered at $x \in X$, i.e., the set $\{y \in X: d(y,x)= r\}$,\\
 $H(A,r)$ & the set $\{x \in X: d(x,A)= r\}$, where $A \subset X$ and $r\geq 0$,
\end{tabular}

\end{center}
\theoremstyle{definition}
\begin{definition}
 The set $U \subset X$ is a \textit{neighborhood} of $M$ if it contains an open set $V \subset X$ which contains $M$. $U$ is an \textit{open neighborhood} of $M$ if it is open and is a neighborhood of $M$.
\end{definition}

\begin{lemma}\thlabel{L1}
Let $O \subset X$ be an open neighborhood of the compact set $M \subset X$. Then there exists an $\epsilon>0$ such that $B(M,\epsilon)=\{x\in X: d(x,M)<\epsilon\}$ is contained in $O$.
\end{lemma}

\begin{proof}
For each $m\in M$, $r_m>0$ is chosen such that $m \in B(m,r_m) \subset V$. Obviously $\{B(m,r_m/2)\}_{m\in M}$ is an open cover of $M$ which is contained in $V$. So by compactness of $M$, there exists a finite subcover $\{B(m_i,r_{m_i}/2)\}_{i=1,2,...,n}$ of $M$. Let $\epsilon=\frac{1}{2}\min_{i=1,2,...,n}{r_{m_i}}$. Then $B(M,\epsilon)$ is contained in $O$; Since if $d(x,M)=\inf_{y\in M}d(x,y)<\epsilon$, $x \in X$, there exists $m\in M$ such that $d(x,m)<\epsilon$, which implies 
$$d(x,m_i)\leq d(x,m)+d(m,m_i)<\epsilon+r_{m_i}/2\leq r_{m_i}$$
for some $i \in {1,2,...,n}$. Hence $x \in B(m_i,r_{m_i})\subset O$.
\end{proof}

The following characterization of neighborhoods of compact sets (whose counterpart for the singletons is trivial) will be handy in developing many results in the future.

\begin{proposition} \thlabel{P2}
The set $U \subset X$ is a neighborhood of the compact set $M \subset X$ if and only if there exists an $\epsilon>0$ such that $B(M,\epsilon)=\{x\in X: d(x,M)<\epsilon\}$ is contained in $U$. 
\end{proposition}

\begin{proof}
($\Leftarrow$) The function $f: X\rightarrow \mathbb{R}$ with $f(x)=d(x,M)$ (sometimes called \textit{the function of distance} to $M$) can easily be shown that is continuous. Then by hypothesis $B(M,\epsilon)=f^{-1}\left((-\infty,\epsilon)\right)$ is an open set contained in $U$. So $U$ is a neighborhood of $M$.\\
($\Rightarrow$) $U$ is a neighborhood of $M$. So it contains an open neighborhood $V$ of $M$. Hence \thref{L1} implies there exists an $\epsilon>0$ such that $B(M, \epsilon) \subset V \subset U$.
\end{proof} 

\theoremstyle{definition}
\begin{definition}
A \textit{dynamical system} on the space $X$ is a triplet $(X,\mathbb{R},\phi)$, where $\phi: X \times \mathbb{R} \rightarrow X$ is a map with the following properties:\\
\begin{align*}
 &(i)\quad \phi(x,0)=x, \quad \text{ for all }x \in X,\\
 &(ii)\quad \phi(\pi(x,t_1),t_2)=\phi(x,t_1+t_2), \quad \text{ for all }t_1,t_2 \in \mathbb{R}\text{ and } x \in X,\\
 &(iii)\quad \phi \text{ is continuous}.
\end{align*}
\end{definition}

We may frequently use the notation $\phi(S,T)$, with $S \subset X$ and $T \subset \mathbb{R}$, to denote the set $\{\phi(s,t): s \in S, t \in T\}$. Also in the case that either $S=\{s\}$ or $T=\{t\}$, we prefer the use of notations $\phi(s,T)$ and $\phi(S,t)$ over $\phi(\{s\},T)$ and $\phi(S,\{t\})$, respectively.

\theoremstyle{definition}
\begin{definition}
The set valued functions $\gamma, \gamma^{+}, \gamma^{-}: X \rightarrow 2^X$ are defined by
\begin{align*}
    \gamma(x)&=\{\phi(x,t): t \in \mathbb{R}\},\\
    \gamma^+{}(x)&=\{\phi(x,t): t \in \mathbb{R}^+\},\\
    \gamma^{-}(x)&=\{\phi(x,t): t \in \mathbb{R}^-\},
\end{align*}
for all $x \in X$. The sets $\gamma(x)$, $\gamma^+(x)$, and $\gamma^-(x)$ are called the \textit{trajectory}, \textit{ positive semi-trajectory}, and \textit{negative semi-trajectory} through $x$, respectively. Furthermore, with the notation introduced earlier, we may readily acknowledge the facts that $\gamma(x)=\phi(x,\mathbb{R})$, $\gamma^+(x)=\phi(x,\mathbb{R}_+)$, and $\gamma^-(x)=\phi(x,\mathbb{R}_{-})$.
\end{definition}

\theoremstyle{definition}
\begin{definition}
 The set valued function $\Omega: X \rightarrow 2^X$ is defined by
\begin{equation*}
\Omega(x)=\{y\in X: \text{there exists a sequence } \{t_n\}_{n\geq 1} \subset \mathbb{R}\text{ such that } t_n \rightarrow +\infty \text{ and } \phi(x,t_n)\rightarrow y\},
\end{equation*}
for any $x \in X$. The set $\Omega(x)$ is called the $\Omega$-limit of $x$.
\end{definition}

\theoremstyle{definition}
\begin{definition}
A set $M \subset X$ is called \textit{invariant under }$\phi$ (or $\phi$-invariant) if 
\begin{equation*}
    \phi(x,t) \in M,
\end{equation*}
for all $x \in M$ and $t \in \mathbb{R}$. It is called \textit{positively invariant under }$\phi$ if the above condition holds for all $t \in \mathbb{R}_{+}$, and is called \textit{negatively invariant under }$\phi$ if the condition holds for all $t \in \mathbb{R}_{-}$. 
\end{definition}
We shall frequently be dealing with transformations $F: X \rightarrow 2^X$. Then for a given subset $M$ of $X$, let $F(M)=\cup \{F(x): x \in M\}=\cup_{x \in M} F(x)$. That having been stated, we can readily have the following statements equivalent to the definitions above:

\begin{proposition} \thlabel{P3}
A set $M \subset X$ is invariant, positively invariant, or negatively invariant under $\phi$ if and only if $\gamma(M)=M$, $\gamma^{+}(M)=M$, or $\gamma^{-}(M)=M$, respectively.
\end{proposition}

\begin{proof}
We prove the proposition only for the invariance property, as the proofs for others are trivial. First, it can be acknowledged that 
$$M \subset \gamma(M)=\cup_{x \in M} \gamma(x),$$
since $x=\phi(x,0) \in \gamma(x)$ for every $x \in M$. Also, $M$ being invariant is equivalent to $\gamma(x) \subset M$ for all $x \in M$, which itself is equivalent to $\gamma(M) \subset M$. This concludes $\gamma(M)=M$.
\end{proof}

\begin{proposition} \thlabel{P4}
For any $x \in X$, $\Omega(x)$ is a closed invariant set and 
$$\overline{\gamma^{+}(x)}=\gamma^{+}(x) \cup \Omega(x).$$
\end{proposition} 

\begin{proof}
Assume $\{y_n\}_{n \geq 1}$ is a sequence in $\Omega(x)$ converging to some $y \in X$. For each $n \geq 1$, there exists a sequence $\{t^n_m\}_{m\geq 1}$, with $t^n_m \rightarrow \infty$, such that $\phi(x,t^n_m) \rightarrow y_n$. So we may choose a sequence $\{m_n\}_{n\geq 1} \subset \mathbb{R}_+$ such that $t^{n}_{m_n} \geq n$ and $d(\phi(x,t^{n}_{m_n}),y_n) <1/n$ for all $n \geq 1$. Then $t^{n}_{m_n} \rightarrow \infty $ as $n$ approaches infinity, and 
\begin{align*}
    d(y,\phi(x,t^{n}_{m_n})) &\leq d(y,y_n)+d(y_n,\phi(x,t^{n}_{m_n}))\\
                             &<d(y,y_n)+1/n,
\end{align*}
which implies $d(y,\phi(x,t^{n}_{m_n})) \rightarrow 0$ by the fact that $y_n \rightarrow y$. This concludes $y \in \Omega(x)$ and hence $\Omega(x)$ is closed. Now let $y \in \Omega(x)$ and $t \in \mathbb{R}$. So there exists a sequence $\{t_n\}_{n \geq 1} \subset \mathbb{R}$, $t_n \rightarrow \infty$, such that $\phi(x,t_n) \rightarrow y$ and then 
\begin{align*}
    \phi(y,t)=&\phi(\lim_{n\rightarrow \infty} \phi(x,t_n),t)\\
             =&\lim_{n\rightarrow \infty} \phi(\phi(x,t_n),t)\\
             =&\lim_{n\rightarrow \infty} \phi(x,t+t_n)
\end{align*}
by continuity of $\phi$. Since $t+t_n \rightarrow \infty$ as $n \rightarrow \infty$, the equalities above imply $\phi(y,t) \in \Omega(x)$. $y$ and $t$ were chosen arbitrarily , therefore $\Omega(x)$ is invariant. We now want to prove the equality in the statement of the proposition. If $y \in \Omega(x)$, then $y$ is a point of closure of $\phi(x,\mathbb{R}_+)=\gamma^{+}(x)$ by the definition of $\Omega(x)$. So $\overline{\gamma^{+}(x)}\supset \gamma^{+}(x) \cup \Omega(x).$ Now let $y \in \overline{\gamma^{+}(x)}$. By the definition of closure of a set, there is a sequence $\{y_n\}_{n \geq 1} \subset \gamma^{+}(x)$ such that $y_n \rightarrow y$. For each $n \geq 1$, there is some $t_n \in \mathbb{R}_+$ such that $y_n=\phi(x,t_n)$. Then either $t_n \rightarrow \infty$, in which case $y \in \Omega(x)$, or there is a bounded subsequence $\{t_{n_m}\}_{m\geq 1}$ and hence a convergent further subsequence $\{t_{n_{m_k}}\}_{k\geq 1}$, where $t_{n_{m_k}} \rightarrow t$ for some $t \in \mathbb{R}_+$. In the latter case, $\phi(x,t_{n_{m_k}}) \rightarrow \phi(x,t)$ by continuity of $\phi$, which implies $y=\phi(x,t) \subset \gamma^+(x)$ by the fact that $\phi(x,t_n) \rightarrow \phi(x,t)$ and uniqueness of the limit. This concludes $\overline{\gamma^{+}(x)}\subset \gamma^{+}(x) \cup \Omega(x)$. The proposition is now proved.
\end{proof}

\begin{lemma} \thlabel{L4'}
Let $X$ be any metric space (not necessarily locally compact) and let $x \in X$. Then for every $t \in \mathbb{R}$, 
$$\Omega(x)=\Omega(\phi(x,t))=\phi(\Omega(x),t).$$
\end{lemma}

\begin{proof}
Our plan for the proof is to show $\Omega(x)\subset \Omega(\phi(x,t))\subset \phi(\Omega(x),t) \subset \Omega(x)$, for all $x \in X$ and $t \in \mathbb{R}$. Let $u \in \Omega(x)$ and $t \in \mathbb{R}$. There exists a sequence $\{t_n\}$, $t_n \rightarrow \infty$, with $\phi(x,t_n) \rightarrow u$. Then for the sequence $t'_n=t_n-t$, we have $t'_n \rightarrow \infty$ and $\phi(\phi(x,t),t'_n)=\phi(x,t_n)$. Hence $u \in \Omega(\phi(x,t))$, which implies $\Omega(x) \subset \Omega(\phi(x,t))$. Now let $y \in \Omega(\phi(x,t))$ where $t \in \mathbb{R}$. So there exists another sequence $\{t_n\}$, $t_n \rightarrow \infty$, with $\phi(\phi(x,t),t_n) \rightarrow y$. Therefor $\phi(x,t_n)=\phi(\phi(\phi(x,t),t_n),-t) \rightarrow \phi(y,-t)$, by continuity of $\phi$. This implies $y \in \phi(\Omega(x),t)$, and hence $\Omega(\phi(x,t))\subset \phi(\Omega(x),t)$. Finally, it is obvious that $\phi(\Omega(x),t) \subset \Omega(x)$ by positive invariance property of $\Omega(x)$ shown in \thref{P4}.
\end{proof}

\begin{proposition} \thlabel{P5}
Let $\{M_i\}_{i \in I}$ be a collection of $\phi$-invariant subsets of $X$, with $I$ an arbitrary index set. Then their union and their intersection are $\phi$-invariant.
\end{proposition}

\begin{proof}
First let $x \in \cup_{i \in I} M_i$. Then $x \in M_i$ for some $i \in I$. So we have $\gamma(x) \subset M_i \subset \cup_{i \in I} M_i$ by the assumption that $M_i$ is invariant under $\phi$. This implies $\cup_{i \in I} M_i$ is invariant under $\phi$. Now let $x \in \cap_{i \in I} M_i$. Then $x \in M_i$ for all $i \in I$. So we have $\gamma(x) \subset \cap_{i \in I} M_i$ by the fact that all $M_i$'s are invariant under $\phi$. This concludes $\cap_{i \in I} M_i$ is invariant under $\phi$.
\end{proof}

\begin{proposition} \thlabel{P5'}
If a set $M \subset X$ is positively invariant, negatively invariant, or invariant, then its closure $\overline{M}$ has the same property.
\end{proposition}

\begin{proof}
We prove the proposition only for invariance property, as the proof for other properties is almost the same. Let $x \in \overline{M}$ and $t \in \mathbb{R}$. Then there exists a sequence $\{x_n\}_{n \geq 1}$ in $M$ such that $x_n \rightarrow x$ and then $\phi(x_n,t) \rightarrow \phi(x,t)$ by continuity of $\phi$. Also, for all $n\geq 1$, $\phi(x_n,t)$ is in $M$ and hence in $\overline{M}$ by the invariance property of $M$. Therefore $\phi(x,t) \in \overline{M}$, as $\overline{M}$ is closed. This implies $\overline{M}$ is invariant, as $x$ and $t$ were chosen arbitrarily. 
\end{proof}

\begin{proposition} \thlabel{P5''}
A set $M \subset X$ is invariant under $\phi$ if and only if its complement, $X\setminus M$, is invariant under $\phi$.
\end{proposition}

\begin{proof}
Assume $M$ is invariant and there are $x \in X\setminus M$ and $t \in \mathbb{R}$ such that $\phi(x,t) \notin X\setminus M$. Then $\phi(x,t) \in M$. This implies $x=\phi(\phi(x,t),-t) \in M$, by the invariance of $M$, which contradicts the assumption that $x \in X\setminus M$. Hence $X\setminus M$ is invariant. Now if $X\setminus M$ is invariant, with a similar argument we may show that $M$ is invariant, since $M=X\setminus (X\setminus M)$. 
\end{proof}
\begin{lemma} \thlabel{L6}
If the space $X$ is locally compact, then for any nonempty compact subset $M$ of $X$ there exists $\epsilon>0$ such that $S[M,\epsilon]=\{x\in X: d(x,M)\leq \epsilon\}$ is compact.
\end{lemma}
		
\begin{proof}
For every point $m\in M$, let $U_m$ and $C_m$ denote an open neighborhood of $m$ and a compact set, respectively, such that $U_m \subset C_m$ ($X$ is locally compact). Also, $r_m$ is chosen such that $B(m,r_m)=\{x\in X : d(x,m)<r_m\} \subset U_m$. $\{B(m,r_m/2)\}_{m\in M}$ covers $M$. Since $M$ is compact, there exists a finite subcover $\{B(m_i,r_{m_i}/2)\}_{i=1,2,...,n}$ of $M$. Choose $\epsilon<\frac{1}{2}\min_{i\in\{1,2,...,n\}}r_{m_i}$. If $x \in S[M,\epsilon]$, i.e. $d(x,M)\leq \epsilon$, since $M$ is compact (and hence sequentially compact) there exists $z \in M$ such that $d(x,z)\leq \epsilon$. $z \in B(m_i,r_{m_i}/2)$ for some $i=1,2,...,n$. So
$$d(x,m_i)\leq d(x,z)+d(z,m_i)<\epsilon+r_{m_i}/2 \leq r_{m_i}.$$
therefore $x \in B(m_i,r_{m_i}) \subset C_{m_i} \subset C$, where $C=\cup_{i\in \{1,2,...,n\}}C_{m_i}$. Hence $S[M,\epsilon] \subset C$. Now $C$ is compact (finite union of compact sets is compact) and $S[M,\epsilon]$ is closed, so $S[M,\epsilon]$ is compact.
\end{proof} 

\begin{lemma} \thlabel{L7}

For any neighborhood $N$ of the nonempty and compact set $M \subset X$, there exists an $\epsilon>0$ such that $S[M,\epsilon]$ is compact and is contained in $N$.
\end{lemma}

\begin{proof}
By \thref{L6}, we can choose $r>0$ such that $S[M,r]$ is compact. Also by \thref{P2}, we can choose $r'>0$ such that $B(M,r') \subset N$. Thus if we pick an $\epsilon>0$, with $\epsilon< \min\{r,r'\}$, then $S[M,\epsilon]$ is first compact, since it is a closed subset of a compact set $S[M,r]$, and second, contained in $B(M,r')$ and hence in $N$.
\end{proof}

\begin{proposition} \thlabel{P8}
In a locally compact metric space $X$, for each $x \in X$, $\Omega(x)$ is nonempty and compact if and only if $\overline{\gamma^{+}(x)}$ is compact.
\end{proposition}

\begin{proof}
($\Leftarrow$) If $\overline{\gamma^{+}(x)}$ is compact then $\Omega(x)$ is compact, since by \thref{P4} it is a closed subset of $\overline{\gamma^{+}(x)}.$  \\
($\Rightarrow$) Assume $\Omega(x)$ is nonempty and compact. Then by \thref{L6} there exists an $\epsilon>0$ such that $S[\Omega(x),\epsilon]$ is compact. Further, there exists a $T\in \mathbb{R}_{+}$ such that $\gamma^{+}(\phi(x,T)) \subset S[\Omega(x),\epsilon]$. For otherwise we can have a sequence $\{t_n\}_{n\geq 1}$ in $\mathbb{R}_{+}$ with $t_n \rightarrow +\infty$ such that $\{\phi(x,t_n)\}_{n
\geq 1} \subset$  $X \setminus S[\Omega(x),\epsilon]$ and, as a result, another sequence $\{t'_n\}_{n\geq 1}$ in $\mathbb{R}_{+}$ with $t'_n \rightarrow +\infty$ such that $\{\phi(x,t'_n)\}_{n
\geq 1} \subset$  $\partial S[\Omega(x),\epsilon]$ due to the fact that $\Omega(x)$ is nonempty. Now since $\partial S[\Omega(x),\epsilon]$ is closed and hence compact, there is a further subsequence $\{t'_{n_k}\}_{k\geq 1}$ of $\{t'_n\}_{n\geq 1}$ for which $\phi(x,t'_{n_k}) \rightarrow y \in$ $\partial S[\Omega(x),\epsilon]$. This contradicts the fact that $y \in \Omega(x)$; since if $y \in$ $\partial S[\Omega(x),\epsilon]$, there is a sequence $\{y_n\}_{n\geq 1}$ in $X \setminus S[\Omega(x),\epsilon]=\{x\in X: d(x,\Omega(x))>\epsilon\}$ such that $y_n \rightarrow y$, which dictates $d(y,\Omega(x))\geq \epsilon$. Also $\gamma^{+}(x)=\phi(x,[0,T]) \cup \gamma^{+}(\phi(x,T))$ where $\phi(x,[0,T]):= \cup_{t \in [0,T]}\{\phi(x,t)\}$ is the flow from $t=0$ to $T$. So $\overline{\gamma^{+}(x)}=\phi(x,[0,T]) \cup \overline{\gamma^{+}(\phi(x,T))}$. Therefore $\overline{\gamma^{+}(x)}$ is compact since $\overline{\gamma^{+}(\phi(x,T))}$, as a closed subset of $S[\Omega(x),\epsilon]$, is compact and so is $\phi(x,[0,T])$ due to continuity of $\phi(x,t)$ in $t$.
\end{proof} 

\begin{lemma} \thlabel{L9}
If $\overline{\gamma^{+}(x)}$ is compact, $x \in X$, then we have 
$$\lim_{t \rightarrow \infty} d(\phi(x,t),\Omega(x))=0.$$ 
\end{lemma}

\begin{proof}
We would like to show this by contradiction. Suppose there is a sequence $\{t_n\}_{n \geq 1}$ in $\mathbb{R}_+$ with $t_n \rightarrow \infty$ such that $d(\phi(x,t_n),\Omega(x)) \nrightarrow 0$. Then there exists an $\epsilon>0$ and a subsequence $\{t_{n_k}\}_{k \geq 1}$ of $\{t_n\}_{n \geq 1}$ such that $d(\phi(x,t_{n_k}),\Omega(x))>\epsilon$ for all $k \geq 1$. Since $\{\phi(x,t_{n_k})\}_{k \geq 1} \subset \gamma^+(x)\subset \overline{\gamma^+(x)}$, and $\overline{\gamma^+(x)}$ is compact, then there can be extracted a further subsequence $\{t_{n_{k_m}}\}_{m \geq 1}$ which converges to some $y \in \Omega(x)$. This contradicts the assumption that $d(\phi(x,t_{n_k}),\Omega(x))>\epsilon$, $k\geq 1$.
\end{proof}

\begin{definition} \thlabel{D6}
For a non-empty compact set $M\subset X$, the sets
\begin{align*}
    &A_\omega(M)=\{x \in X: \Omega(x) \cap M \neq \emptyset\}\\
    &A(M)=\{x \in X: \Omega(x)\neq \emptyset \text{ and } \Omega(x) \subset M\}
\end{align*}
are respectively called the \textit{region of weak attraction}, and \textit{the region of attraction} of $M$. A point $x$ in $A_{\omega}(M)$ or $A(M)$ is said to be \textit{weakly attracted to} $M$ or \textit{attracted to} $M$, respectively.
\end{definition} 

\begin{proposition} \thlabel{P10}
For any set $M \subset X$, $A(M)$ is invariant.
\end{proposition}

\begin{proof}
Let $x \in A(M)$. $\Omega(\phi(x,t))= \Omega(x)$ for all $t \in \mathbb{R}$ by \thref{L4'}. This implies $\Omega(\phi(x,t)) \neq \emptyset$ and $\Omega(\phi(x,t)) \subset M$ by \thref{D6}, and hence $\Omega(x,t) \in A(M)$ for all $t \in \mathbb{R}$. Therefore $A(M)$ is invariant under $\phi$. 
\end{proof}

\begin{definition} \thlabel{D5}
A given set $M$ is said to be \textit{weakly attracting} (or \textit{a weak attractor}) if $A_{\omega}(M)$ is a neighborhood of $M$, and \textit{attracting}(or an \textit{attractor}) if $A(M)$ is a neighborhood of $M$.
\end{definition}

\begin{proposition} \thlabel{P11}
Given $M$, a point $x \in X$ is weakly attracted to $M$ if and only if there is a sequence $\{t_n\}_{n \geq 1}$ in $\mathbb{R}_+$ with $t_n \rightarrow \infty$ such that $d(\phi(x,t_n),M)\rightarrow 0$.
\end{proposition}

\begin{proof}
($\Rightarrow$) Suppose for a point $x \in X$, $\Omega(x)\cap M\neq \emptyset$. Then there exists $y \in M$ and a sequence $\{t_n\}_{n\geq 1}$ in $\mathbb{R}_+$ such that $t_n \rightarrow \infty$ and  $d(\phi(x,t_n),y)\rightarrow 0$, hence $d(\phi(x,t_n),M)\rightarrow 0$.\\
($\Leftarrow$) Now suppose there is a sequence $\{t_n\}_{n\geq 1}$ in $\mathbb{R}_+$ with $t_n \rightarrow \infty$ such that $d(\phi(x,t_n),M)\rightarrow 0$. So for $\epsilon=\frac{1}{m}>0$, $m \in \mathbb{N}$, there exists an $N\geq 1$ such that $d(\phi(x,t_n),M)=\inf_{y \in M} d(\phi(x,t_n),y)<\frac{1}{m}$ for $n\geq N$. This implies we can choose a sequence $\{y_m\}_{m \geq 1}$ in $M$ and a subsequence $\{t_{n_m}\}_{m\geq 1}$ of $\{t_n\}_{n\geq 1}$ such that $d(\phi(x,t_{n_m}),y_m)\rightarrow 0$ as $m\rightarrow \infty$. Since $M$ is compact (and hence sequentially compact) $\{y_m\}_{m\geq 1}$ possesses a subsequence $\{y_{m_k}\}_{k\geq 1}$ convergent to a point $y\in M$. Therefore there is a sequence $\{t'_k\}_{k \geq 1}$ in $\mathbb{R}_+$ with $t'_k=t_{n_{m_k}}$ and $t'_k \rightarrow \infty$ such that $d(\phi(x,t'_k),y)\rightarrow 0$, which implies $\Omega(x)\cap M\neq \emptyset$.
\end{proof}

\begin{proposition} \thlabel{P12}
Given $M \subset X$, a point $x \in X$ is attracted to $M$ if and only if 
$$\lim_{t \rightarrow \infty} d(\phi(x,t),M)=0.$$
\end{proposition}

\begin{proof}
($\Rightarrow$) Assume $x$ is attracted to $M$, i.e., $\Omega(x)\neq \emptyset$ and $\Omega(x)\subset M$. So $\Omega(x)$ is compact since it is closed and $M$ is compact. Then  $\overline{\gamma^{+}(x)}$ is compact by \thref{P8} and thus, $d(\phi(x,t),\Omega(x)) \rightarrow 0$ as $t \rightarrow \infty$ by \thref{L9}. Now since $\Omega(x) \subset M$, we have $d(\phi(x,t),M) \leq d(\phi(x,t),\Omega(x))$. Hence $d(\phi(x,t),M) \rightarrow 0$ as $t \rightarrow \infty$.\\
($\Leftarrow$) If $d(\phi(x,t),M) \rightarrow 0$ as $t \rightarrow \infty$, we can show $\Omega(x) \neq \emptyset$ in a similar fashion to the ($\Leftarrow$) part of the proof of \thref{P11}. In order to show $\Omega(x)\subset M$, first let $y \in \Omega(x)$. So there is a sequence $\{t_n\}_{n\geq 1}$ in $\mathbb{R}_+$ such that $t_n \rightarrow \infty$ and $\phi(x,t_n) \rightarrow y$ (or equivalently $d(\phi(x,t_n),y) \rightarrow 0$). On the other hand $d(\phi(x,t_n),M)\rightarrow 0$ by hypothesis. So, again in a similar fashion to the ($\Leftarrow$) part of the proof of \thref{P11}, we can choose a subsequence $\{t_{n_k}\}_{k \geq 1}$ of $\{t_n\}_{n\geq 1}$ such that $d(\phi(x,t_{n_k}),y')\rightarrow 0$ for some $y' \in M$. $y'$ has no other way but $y'=y$ (since $X$ is  metric, hence Hausdorff.) Therefore $y \in M$. This concludes $\Omega(x)\subset M$.
\end{proof}

\begin{proposition} \thlabel{P13}
If $M$ is an attractor, then its region of attraction, $A(M)$, is open.
\end{proposition}

\begin{proof}
Suppose there exists a point $x \in \partial A(M) \cap A(M)$, where $\partial A(M)$ is the boundary of $A(M)$. Since $M$ is an attractor, $A(M)$ is a neighborhood of $M$, and hence $\partial A(M) \cap M=\emptyset$. Furthermore, \thref{P5'} together with \thref{P10} imply $\overline{A(M)}$ is invariant. Hence, the exterior of $A(M)$, ext$A(M)$, is invariant by the fact that ext$A(M)=X\setminus\overline{A(M)}$ and \thref{P5''}. Now it can be verified that $\partial A(M)$ is closed and positively invariant since 
$$\partial A(M)=\overline{\text{ext}A(M)} \cap \overline{A(M)},$$
and the intersection of two invariant sets is invariant (\thref{P5}). So we have 
$$\Omega(x) \subset \overline{\gamma^+(x)} \subset \partial A(M)$$
by \thref{P4}. But $\Omega(x) \neq \emptyset$ and  $\Omega(x) \subset M$ by the assumption that $x \in A(M)$. This contradicts the fact that $\partial A(M) \cap M=\emptyset$. Therefore $\partial A(M) \cap A(M) =\emptyset$, i.e., $A(M)$ is open.
\end{proof}

\begin{definition}
Let $M\subset X$ be a compact set. $M$ is said to be \textit{Lyapunov stable} under $\phi$ if for every neighborhood $N$ of $M$ there exists a neighborhood $N'$ of $M$ such that $\gamma^+(N')\subset N$, where $\gamma^+(N')$ is defined with respect to $\phi$.
\end{definition}

It should be noted that there are alternative definitions of Lyapunov stability in the literature which turn out to be equivalent to the one stated above, at least in the topology settings we have adopted. We will provide just two of them in the following propositions.

\begin{proposition} \thlabel{P14}
A compact set $M\subset X$ is Lyapunov stable under $\phi$ if and only if every neighborhood $U$ of $M$ contains a positively invariant (under $\phi$) neighborhood $V$ of $M$.
\end{proposition}

\begin{proof}
($\Rightarrow$) Assume $M$ is Lyapunov stable. Then for a given neighborhood $U$ of $M$ there exists a neighborhood $U'$ of $M$ such that $\gamma^+(U') \subset U$. Now $V=\gamma^+(U')$ is a positively invariant neighborhood of $M$  contained in $U$.\\
$(\Leftarrow)$ Suppose $N$ is a neighborhood of $M$. $N$ contains a positively invariant neighborhood $N'$ of $M$, so $\gamma^+(N')=N' \subset N$. Hence $M$ is Lyapunov stable.
\end{proof}

\begin{proposition} \thlabel{P15}
A compact set $M\subset X$ is Lyapunov stable if and only if for every $\epsilon>0$ there exists $\delta=\delta(\epsilon)>0$ such that if $d(x,M)<\delta$, $x \in X$, then $d(\phi(x,t),M)<\epsilon$ for all $t\geq 0$; or equivalently 
$$\gamma^+(B(M,\delta)) \subset B(M,\epsilon).$$
\end{proposition} 

\begin{proof} ($\Rightarrow$) Assume $M$ is stable. Then for every neighborhood $B(M,\epsilon)$ of $M$ there exists a neighborhood $N$ of $M$ which its contains, by \thref{P2}, an open neighborhood $B(M, \delta)$ of $M$ such that $B(M, \delta) \subset \gamma^+(N) \subset B(M,\epsilon)$. Thus $\gamma^+(B(M,\delta)) \subset \gamma^+(N) \subset B(M,\epsilon)$. \\
($\Leftarrow$) Let $U$ be a neighborhood of $M$. Then there exists an $\epsilon>0$ such that $B(M,\epsilon) \subset U$, again by \thref{P2}. Let $\delta>0$ respond to the $\epsilon$ challenge stated above. Clearly $B(M,\delta)$ is a neighborhood of $M$ for which $\gamma^+(B(M,\delta)) \subset B(M,\epsilon) \subset U$. $U$ was chosen arbitrarily, hence $M$ is stable.
\end{proof}

The following Proposition reveals the connection between Lyapunov stability and positive invariance of compact subsets of $X$.\\

\begin{proposition} \thlabel{P16}
If a set $M \subset X$ is Lyapunov stable, then it is positively invariant.
\end{proposition}

\begin{proof} 
If $M$ is stable, it is the intersection of all of its positively invariant neighborhoods. To prove so, assume $x \in X$ is such that it belongs to all positively invariant neighborhoods of $M$ but not to $M$. Take $r=d(x,M)>0$. Since $M$ is stable, there exists a positively invariant neighborhood $N$ of $M$, contained in $B(M,r/2)$. By construction, $x \notin B(M,r/2)$, thus $x \notin N$; this contradicts the assumption that $x$ belongs to all positively invariant neighborhoods of $M$. Now, \thref{P5} implies that $M$ is positively invariant.
\end{proof}

\begin{definition}
A compact set $M \subset X$ is \textit{(Lyapunov) asymptotically stable} if it is (Lyapunov) stable and attracting.
\end{definition}

\begin{proposition} \thlabel{P17}
If a set $M \subset X$ is (Lyapunov) asymptotically stable under $\phi$, with region of attraction $A(M)$, then for every compact set $K$ in $A(M)$ and every neighborhood $N$ of $M$ there exists $T>0$ such that $\phi(K,[T,+\infty)) \subset N$.
\end{proposition}

\begin{proof}
The property mentioned in the proposition, with $A(M)$ being replaced with an arbitrary set $U \subset X$, is equivalent to a property usually called, in the literature, \textit{uniform attraction relative to }$U$. So the proposition above can be translated to the following: If a set $M \subset X$ is (Lyapunov) asymptotically stable, with region of attraction $A(M)$, then it is uniformly attracting relative to $A(M)$. Look at reference 1, Theorem 5.5 in Chapter 5, for a complete proof. 
\end{proof}

Note that when $U$ is a neighborhood of $M$, which is the case for $U=A(M)$ when $M$ is asymptotically stable, relative attraction relative to $U$ reduces to another property named \textit{uniform attraction}. For an extensive exposition see reference 1, Definition 5.1 in Chapter 5.  \\

\begin{theorem} \thlabel{T18}
A compact set $M \subset X$ is (Lyapunov) asymptotically stable if and only if there exists a continuous real valued function $L$ defined on a neighborhood $N$ of $M$ such that\\

\noindent(i) $L(x)=0$ if $x \in M$, and $L(x)>0$ if $x \notin M$;\\
(ii) $L(\phi(x,t))<L(x),$ for all $ x \notin M$ and $t>0$ with $\phi(x,[0,t]) \subset N$.
\end{theorem}

Before proving the theorem we need to provide a few Lemmas, one of which is a classic topological fact on sets in general, another one on invariance of a set $M$ as described in \thref{T18}, and the rest are facts from calculus.\\ 

\begin{lemma} \thlabel{L19}
Every path $s: [0,1] \rightarrow X$ from a point in $A$, $A \subset X$, to a point not in $A$ passes through the boundary $\partial A$ of $A$, i.e., $$s([0,1]) \cap \partial A \neq \emptyset.$$
\end{lemma}

\begin{proof}
If either $s(0) \in \partial A$ or $s(1) \in \partial A$, then $s$ is passing through $\partial A$. Now assume that $s(0) \in A^{\mathrm{o}}$ and $s(1) \in X\setminus \overline{A}$, where $A^{\mathrm{o}}$ is the interior of $A$ and $\overline{A}$ is the closure of $A$, and that $s([0,1]) \cap \partial A=\emptyset$. Since $A^{\mathrm{o}}$ and $X\setminus \overline{A}$ are open and the path $s$ is continuous, we have $s^{-1}(A^{\mathrm{o}})=U \cap [0,1]$ and $s^{-1}(X\setminus \overline{A})=V \cap [0,1]$ for some open sets $U$ and $V$ in $\mathbb{R}$. Then, the sets $U$ and $V$ separate $[0,1]$. This follows from the facts that $A^{\mathrm{o}}$ and $X\setminus \overline{A}$ are disjoint, $[0,1]=s^{-1}(A^{\mathrm{o}}) \cup s^{-1}(X\setminus \overline{A}) \cup s^{-1}(\partial A)$, and $s^{-1}(\partial A)=\emptyset$ by the assumption made earlier. The fact that $[0,1]$ is separated by $U$ and $V$ contradicts  that $[0,1]$ is connected. Therefore, $s([0,1]) \cap \partial A \neq \emptyset.$
\end{proof}

\begin{lemma} \thlabel{L20}
A compact set $M \subset X$, with the the properties (i) and (ii) in the sufficiency part of \thref{T18}, is positively invariant under $\phi$.
\end{lemma}

\begin{proof}
Let $x \in M$ and $t>0$ such that $\phi(x,t) \notin M$. Then $L(\phi(x,t))>0$ by hypothesis (i) in \thref{T18}. Since $L$ and $\phi$ are continuous, Intermediate Value Theorem implies there exists a $t' \in (0,t)$ such that $0<L(\phi(x,t'))<L(\phi(x,t))$. So $\phi(x,t') \notin M$ again by hypothesis (i), and hence $L(\phi(x,t))=L(\phi(\phi(x,t'),t-t'))< L(\phi(x,t'))$ by hypothesis (ii) of the same theorem. This contradiction concludes that $M$ is positively invariant under $\phi$.
\end{proof}

\begin{lemma} \thlabel{L21} 
Let $K \subset X$ where $X$ is an arbitrary metric space. Assume $f$ is a continuous real valued function on $K$ such that $f(\phi(x,t))\leq f(x)$ for all $x \in K$ and $t \geq 0$ with $\phi(x,[0,t]) \subset K$. Then if for some $x$, $\overline{\gamma^+(x)} \subset K$, we have $f(y)=f(z)$ for every $y,z \in \Omega(x)$. 
\end{lemma}

\begin{proof}
Suppose there exist $y,z \in \Omega(x)$ such that $f(y)<f(z)$. Then there are sequences $\{t_n\}_{n \geq 1}$ and $\{t'_n\}_{n \geq 1}$ in $\mathbb{R}_+$, with $t_n \rightarrow \infty $ and $t'_n \rightarrow \infty $, such that $\phi(x,t_n) \rightarrow y$ and $\phi(x,t'_n) \rightarrow z$. We may extract a subsequence $\{t'_{n_m}\}_{m \geq 1}$ of $\{t'_n\}_{n \geq 1}$ such that $t'_{n_m} \geq t_m$ for all $m \geq 1$. Then we have
\begin{align*}
    f(\phi(x,t'_{n_m})) &= f(\phi(\phi(x,t_m),t'_{n_m}-t_m))\\
                        &\leq f(\phi(x,t_m)),
\end{align*}
since $t'_{n_m}-t_m \geq 0$, and $\phi(\phi(x,t_m),[0,t'_{n_m}-t_m]) \subset \overline{\gamma^+(x)} \subset K$. This leads to $f(z) \leq f(y)$ by continuity of $f$, which contradicts the initial assumption that $f(y) < f(z)$. So the Lemma is proved.
\end{proof}

\begin{lemma} \thlabel{L22} 
If $f,g: \mathcal{A} \rightarrow \mathbb{R}$ are bounded functions, where $\mathcal{A}$ is an arbitrary set, then 
$$|\sup_{x \in \mathcal{A}} f(x)-\sup_{x \in \mathcal{A}} g(x)| \leq \sup_{x \in \mathcal{A}} |f(x)-g(x)|.$$
\end{lemma}

\begin{proof}
we know that for all $x \in \mathcal{A}$, 
\begin{align}
    f(x)&=g(x)+(f(x)-g(x)) \notag\\
    \Rightarrow \sup_{x \in \mathcal{A}} f(x)&=\sup_{x \in \mathcal{A}}\left[g(x)+(f(x)-g(x))\right] \notag\\
        &\leq \sup_{x \in \mathcal{A}}g(x) + \sup_{x \in \mathcal{A}} (f(x)-g(x)) \notag\\
        &\leq \sup_{x \in \mathcal{A}} g(x)+\sup_{x \in \mathcal{A}} |f(x)-g(x)| \notag\\
    \Rightarrow \sup_{x \in \mathcal{A}} f(x)-\sup_{x \in \mathcal{A}} g(x)&\leq \sup_{x \in \mathcal{A}} |f(x)-g(x)| \notag.                  
\end{align}
Likewise, we may prove 
\begin{equation}
    \sup_{x \in \mathcal{A}} g(x)-\sup_{x \in \mathcal{A}} f(x)\leq \sup_{x \in \mathcal{A}} |g(x)-f(x)|, \notag
\end{equation}
which concludes 
$$|\sup_{x \in \mathcal{A}} f(x)-\sup_{x \in \mathcal{A}} g(x)| \leq \sup_{x \in \mathcal{A}} |f(x)-g(x)|.$$
\end{proof}

\begin{proposition} \thlabel{P23}
Let $f: A\times B \rightarrow \mathbb{R}$ be a continuous function, where $A$ and $B$ are compact subsets of metric spaces $(X,d)$ and $(Y,d')$ respectively. Then the function $g:A \rightarrow \mathbb{R}$, with $$g(x) = \sup_{y \in B} f(x,y)$$ is uniformly continuous.
\end{proposition}

\begin{proof}
$A$ and $B$ are compact, so $f$ is uniformly continuous. Let $\delta>0$ respond to the $\epsilon/2$ challenge for uniform continuity of $f$. Then, if $d(x,x')<\delta$, we have
\begin{align}
    |g(x')-g(x)|&=|\sup_{y \in B} f(x',y)-\sup_{y \in B} f(x,y)| \notag\\
                &\leq \sup_{y \in B} |f(x',y)-f(x,y)| \label{eq1}\\
                &\leq \epsilon/2 \notag\\
                &<\epsilon, \notag
\end{align}
where inequality (\ref{eq1}) is valid by \thref{L22}; since continuity of $f$ together with compactness of $B$ imply the function $y \rightarrow f(x,y)$ to be bounded for every $x \in A$. Thus $g$ is uniformly continuous.
\end{proof}

\noindent\textit{Proof of \thref{T18}.}
 ($\Leftarrow$) In the light of \thref{P14}, in order to establish stability we will show that every neighborhood $U$ of $M$ contains a positively invariant neighborhood $V$ of $M$. Let $U'=U\cap N$. Trivially, $U'$ is a neighborhood of $M$. So by \thref{L7}, we can choose $\epsilon>0$ such that $S[M,\epsilon] \subset U'$ and $S[M,\epsilon]$ is compact. Pick $\epsilon' \in (0,\epsilon]$ such that $H(M,\epsilon')=\{x\in X:d(x,M)=\epsilon'\}$ is non-empty. Let $\alpha=\min \{L(x): x \in H(M,\epsilon')\}$. Obviously $L(x)>0$ for all $x \in H(M,\epsilon')$, so $\alpha>0$ by continuity of $L$ and compactness of $H(M,\epsilon')$ (Note that $H(M,\epsilon') \subset S[M,\epsilon]$.) Take $\beta \in (0,\alpha)$ and let $$V=\{x \in S[M,\epsilon']: L(x)\leq \beta\}.$$
$M$ is positively invariant by \thref{L20}. Now if $x \in V\setminus M$, then $L(\phi(x,t))<L(x)\leq \beta$ for all $t>0$ by hypothesis (ii), and also $d(\phi(x,t),M)$ can not be greater than $\epsilon'$ because if it were the case for some $t>0$, then there would exist a $t' \in [0,t)$ such that $\phi(x,t') \in$ $\partial S[M,\epsilon'] \subset H(M,\epsilon')$ by \thref{L19}, which would cause $L(\phi(x,t'))\geq \alpha>\beta\geq L(x)$. This implies for every $x \in V\setminus M$, $\phi(x,t)$ lies in $V$ for all $t>0$, which together with positive invariance of $M$ conclude that $V$ is positively invariant. On the other hand, $M \subset \{x \in B(M,\epsilon'): L(x)< \beta\} \subset V$. So $V$ is a positively invariant  neighborhood of $M$ contained in $U$. This establishes that $M$ is stable. Now it remains to show $M$ is an attractor. Since $S[M,\epsilon]$ is compact, the set $V$ constructed above is a compact positively invariant neighborhood of $M$. Hence for any $x \in V$, first $\Omega(x) \neq \emptyset$, by the fact that compactness implies sequential compactness; and second, $L$ is constant on $\Omega(x)$ by \thref{L21}, thus $\Omega(x) \subset M$; since if there is $y \in \Omega(x)$ such that $y \notin M$, then for any $t>0$, $\phi(y,t) \in \Omega(x)$ by \thref{P4}, and $L(\phi(y,t))<L(y)$, by hypothesis (ii), which are in contradiction with the fact that $L$ is constant on $\Omega(x)$. Hence $V \subset A(M)$. This concludes, in this case, that $M$ is asymptotically stable. For the case when $H(M,\epsilon')=\emptyset$ for all $\epsilon' \in (0,\epsilon]$, $M$ will be a compact positively invariant neighborhood of itself by \thref{L20}, which readily implies $M$ is stable and, with a similar argument to the last case, $M \subset A(M)$. Therefore  $M$ is asymptotically stable. \\
($\Rightarrow$) Now suppose $M$ is asymptotically stable and $A(M)$ is its region of attraction. Let us define $\ell: A(M) \rightarrow \mathbb{R}_+$ with
$$\ell(x):=\sup_{t \in [0, +\infty)}\{d(\phi(x,t),M)\}.$$
Note that the function $\ell$ is well defined: $\ell(x)=d(x,M)=0$, for $x \in M$,  since $M$ is stable and hence positively invariant(\thref{P16}). For the case $x \notin M$, let $d_x=d(x,M)>0$. By asymptotic stability of $M$, there exists $T_x \geq 0$ such that 
\begin{equation} \label{eq2}
    d(\phi(x,t),M)<d_x 
\end{equation}
for all $t> T_x$, and then we have
\begin{align*}
    \ell(x)&=\max \left\{\sup_{t \in [0,T_x]}\{d(\phi(x,t),M)\},\sup_{t \in (T_x,\infty)}\{d(\phi(x,t),M)\}\right\}\\
           &=\sup_{t \in [0,T_x]}\{d(\phi(x,t),M)\}\\
           &<\infty
\end{align*}
as a result of \eqref{eq2}, the fact that 
$$\sup_{t \in [0,T_x]}\{d(\phi(x,t),M)\} \geq d(x,M)=d_x,$$
and of continuity of $\phi(x,t)$ in $t$. $A(M)$ is invariant (\thref{P10}), so $\ell(\phi(x,t))$ is well defined for all $t\geq 0$ and $x \in A(M)$. Also, $\ell(\phi(x,t))$ is non-increasing in $t \in [0,\infty)$: for every $t_2>t_1\geq 0$ and $x \in A(M)$ we have
\begin{align*}
\ell(\phi(x,t_2))=&\sup_{t \in [0, +\infty)}\{d(\phi(\phi(x,t_2),t),M)\}\\
            =&\sup_{t \in [0,+\infty)}\{d(\phi(x,t_2+t),M)\}\\
            =&\sup_{t \in [t_2,+\infty)}\{d(\phi(x,t),M)\}\\
         \leq& \sup_{t \in [t_1, +\infty)}\{d(\phi(x,t),M)\}\\
            =&\ell(\phi(x,t_1)).
\end{align*}
So as a special case, $\ell(\phi(x,t))\geq \ell(\phi(x,0))=\ell(x)$ for $t>0$. Also we claim that $\ell$ is continuous, first on $M$, and then on $A(M)\setminus M$. Take $x \in M$. Corresponding to the stability criterion in \thref{P15}, assume $\delta >0$ responds to the $\epsilon/2$ challenge for stability of $M$. So $d(y,x)<\delta$ implies $d(y,M)<\delta$ and then $d(\phi(y,t),M)< \epsilon/2$, for $t \geq 0$. Therefore we have 
\begin{align*}
                |\ell(y)-\ell(x)|=\ell(y)&=\sup_{t\geq 0}\{d(\phi(y,t),M)\}\\
                                &\leq \epsilon/2\\
                                &<\epsilon,
\end{align*}
which implies $\ell$ is continuous on $M$. Now let $x \notin M$ and $d_x=d(x,M)$. Since $X$ is locally compact, there are an $r_x>0$, a compact neighborhood $C_x$ of $x$, and an open neighborhood $U_x$ of $x$ such that $x \in B(x,r_x) \subset U_x \subset C_c$. Also since $A(M)$ is open, we can choose $r'_x,r''_x$, with $0<r''_x<r'_x<\min \{r_x,d_x/2\}$, so that $S[x,r''_x]$ is compact due to $S[x,r''_x]\subset B(x,r_x) \subset C_x$,  and that $S[x,r''_x] \subset B(x,r'_x) \subset A(M)$. Hence \thref{P17} implies that there is a $T>0$ such that 
\begin{equation} \label{eq3}
    \phi(S[x,r''_x],[T,\infty)) \subset B(M,d_x/2).
\end{equation} 
Now if $y \in B(x,r''_x)$, we have 
\begin{equation} \label{eq4}
    \sup_{0 \leq t < \infty}d(\phi(y,t),M) \geq d(y,M)>d_x/2,
\end{equation}
since
\begin{align*}
    d(y,M)&=\inf_{m \in M}d(y,m) \\
          &\geq\inf_{m \in M} \{d(x,m)-d(y,x)\} \\
          &=\inf_{m \in M} \{d(x,m)\}-d(y,x) \\
          &=d(x,M)-d(y,x) \\
          &>d_x-r''_x \\
          &>d_x-d_x/2 \\
          &=d_x/2.
\end{align*}
This leads to
\begin{align}
    |\ell(y)-\ell(x)|=&|\sup_{0 \leq t < \infty}d(\phi(y,t),M)-\sup_{0 \leq t <\infty}d(\phi(x,t),M)| \notag\\
                =&|\sup_{0 \leq t \leq T}d(\phi(y,t),M)-\sup_{0 \leq t \leq T}d(\phi(x,t),M)| \label{eq5}\\
             \leq&\sup_{0 \leq t \leq T}|d(\phi(y,t),M)-d(\phi(x,t),M)| \label{eq6}\\
             \leq& \sup_{0 \leq t \leq T}d(\phi(y,t),\phi(x,t)), \label{eq7}
\end{align}
where \eqref{eq5} is true because of the fact that 
$$\sup_{0 \leq t < \infty}d(\phi(y,t),M)=\sup_{0 \leq t < T}d(\phi(y,t),M)$$
by \eqref{eq3} and \eqref{eq4}, and of the fact that
$$\sup_{0 \leq t < \infty}d(\phi(y,t),M) \geq \sup_{0 \leq t \leq T}d(\phi(y,t),M)\geq \sup_{0 \leq t < T}d(\phi(y,t),M);$$
Also \eqref{eq6} is true since for every $y \in X$, the function $t \rightarrow d(\phi(y,t),M)$ is continuous and hence bounded on $[0,T]$. Again since $\phi$ is continuous, the function $f: S[x,r''_x] \times [0,T] \rightarrow \mathbb{R}_+$, with $f(y,t)=d(\phi(y,t),\phi(x,t))$ is continuous. Hence the function 
$$g(y)=\sup_{0 \leq t \leq T}f(y,t)$$
is uniformly continuous by \thref{P23}. So let $\delta>0$ respond to the $\epsilon$ challenge for uniform continuity of $g$. Then we claim $\delta'=\min \{\delta,r''_x\}$ will respond to the $\epsilon$ challenge for continuity of $\ell$ at $x$: if $y \in B(x,\delta')$, inequality \eqref{eq7} implies 
\begin{align*}
    |\ell(y)-\ell(x)|&\leq \sup_{0 \leq t \leq T}d(\phi(y,t),\phi(x,t))\\
                     &=g(y)\\
                     &=|g(y)-g(x)|\\
                     &<\epsilon.
\end{align*}
This concludes $\ell$ is continuous on $A(M)\setminus M$, and consequently on $A(M)$. The only issue with the function $\ell$ is that it may be not strictly decreasing along parts of the trajectories in $A(M)$ that do not lie in $M$, as required in the theorem. Based on the function $\ell$, we can define a function $L: A(M) \rightarrow \mathbb{R}_+$, with 
$$L(x):= \int_0^{\infty} \alpha(t)\ell(\phi(x,t))dt,$$
where $\alpha$ is a positive non-increasing Lebesgue integrable function. Let $x \in A(M)$. So we have  
\begin{equation*}
    \alpha(t)\ell(\phi(x,t)) \leq \alpha(t)\ell(x) 
\end{equation*}
for all $t \geq 0$. This with the assumption that $\alpha$ is integrable imply $L(x)$ is well defined. Furthermore, by \thref{L6} and the fact that $A(M)$ is open, there is a $\delta>0$ such that $S[x,\delta] \subset A(M)$ and $S[x,\delta]$ is compact. Now let $\overline{\ell}:=\max_{y\in S[x,\delta]} \ell(y)<\infty$. Then
\begin{align}
    \lim_{y \rightarrow x} L(y)&=\lim_{y \rightarrow x}\int_0^{\infty} \alpha(t)\ell(\phi(y,t))dt \notag\\
                &=\int_0^{\infty} \lim_{y \rightarrow x} \alpha(t)\ell(\phi(y,t))dt \label{eq8}\\
                &=\int_0^{\infty} \alpha(t)\ell(\phi(x,t))dt\notag\\
                &=L(x). \notag
\end{align}
The change of integration and limit operations in \eqref{eq8} needs to be verified: for any sequence $x_n \rightarrow x$, the sequence $\{\alpha(t)\ell(\phi(x_n,t))\}_{n\geq 1}$ is convergent to $\alpha(t)\ell(\phi(x,t))$ pointwise in $t$, by continuity of $\ell$ and $\phi$; and (with no loss of generality) by assuming $\{x_n\}_{n \geq 1} \subset S[x,\delta]$ we have
$$\alpha(t)\ell(\phi(x_n,t)) \leq \alpha(t)\ell(x_n) \leq \overline{\ell}\alpha(t) \quad\quad \forall n\geq 1.$$
So Lebesgue Dominated Convergence Theorem implies $$\lim_{n \rightarrow \infty} \int_{0}^{\infty} \alpha(t)\ell(\phi(x_n,t))dt=\int_{0}^{\infty} \lim_{n \rightarrow \infty}\alpha(t)\ell(\phi(x_n,t))dt,$$ and hence \eqref{eq8}. Therefore $L$ is continuous. If $x \in M$, $L(x)=0$ since $\ell(\phi(x,t))=0$ for all $t\geq 0$ by monotonicity of $\ell(\phi(x,t))$ in $t$; and if $x \notin M$, due to continuity of $\ell(\phi(x,t))$ in $t$, there is a $\Delta t>0$ such that $\ell(\phi(x,t))>0$ for $t \in [0, \Delta t]$, which leads to 
\begin{align}
    L(x)&=\int_0^{\infty} \alpha(t)\ell(\phi(x,t))dt \notag\\
        &\geq \int_0^{\Delta t} \alpha(t)\ell(\phi(x,t))dt \notag\\
        &>0. \label{eq9}
\end{align}
Note that inequality \eqref{eq9} is strict, since if $\int_0^{\Delta t} \alpha(t)\ell(\phi(x,t))dt=0$, it implies $\alpha(t)\ell(\phi(x,t))=0$ almost everywhere in $t$, which contradicts the strict positivity of the integrand on $[0,\Delta t]$.
Finally, in order to show property (ii), take an $x \notin M$ and a $t>0$. It is clear that $L(\phi(x,t)) \leq L(x)$ by the fact that $\ell(\phi(x,t))$ is non-increasing in $t$. If it were the case that $L(\phi(x,t))=L(x)$, then 
\begin{align}
   \int_0^{\infty} \alpha(t')\ell(\phi(x,t+t'))dt-\int_0^{\infty} \alpha(t')\ell(\phi(x,t'))dt'=0 \notag\\
   \Rightarrow \int_0^{\infty} \alpha(t')[\ell(\phi(x,t+t'))-\ell(\phi(x,t'))]dt'=0 \label{eq10}\\
   \substack{\alpha>0\\\Longrightarrow} \quad  \ell(\phi(x,t+t'))=\ell(\phi(x,t')) \quad\quad\forall t' \in [0,\infty). \label{eq11}
\end{align}
The equality in \eqref{eq11} is valid everywhere in $t$ (not just almost everywhere, as it should be the conclusion for Lebesgue integral), since if $\ell(\phi(x,t+t'))<\ell(\phi(x,t'))$ for some $t'\geq 0$, with a similar argument to that for positivity of $L$ on $A(M)\setminus M$, it would require the equality sign in \eqref{eq10} to turn into $<$ sign. Further, \eqref{eq11} implies that $\ell(\phi(x,nt))=\ell(x)$ for all $n\geq 1$, and hence 
\begin{equation} \label{eq12}
    \lim_{n \rightarrow \infty}\ell(\phi(x,nt))=\ell(x)>0
\end{equation}
when $x \notin M$. However, since $x\in A(M)$, $d(\phi(x,t),M) \rightarrow 0$ as $t \rightarrow \infty$. So for any $\epsilon>0$ there is a $T>0$ such that $d(\phi(x,t),M)<\epsilon/2$, when $t\geq T$. Then we have
\begin{align*}
    \ell(\phi(x,t))&=\sup_{t' \in [0, +\infty)}\{d(\phi(x,t+t'),M)\}\\
                   &=\sup_{t' \in [t, +\infty)}\{d(\phi(x,t'),M)\}\\
                   &\leq \sup_{t' \in [T, +\infty)}\{d(\phi(x,t'),M)\}\\
                   &\leq \epsilon/2\\
                   &<\epsilon
\end{align*}
for $t\geq T$, which implies $\lim_{t \rightarrow \infty}\ell(\phi(x,t))=0$. This contradicts \eqref{eq12}. Therefore $L(\phi(x,t)) < L(x)$ for all $x \notin M$, and $t>0$. \hfill $\qed$
\\

Now, consider the system of ordinary differential equations
\begin{equation} \label{eq13}
    \dot{\mathbf{x}}=V(\mathbf{x}),
\end{equation}
where $V:D \rightarrow \mathbb{R}^n$, with $D$ a domain (an open connected set) in the space $\mathbb{R}^n$, is continuous. We may describe a solution of \eqref{eq13} by a function $\varphi (t,x_0)$ not only of $t \in \mathbb{R}_+$, but also of coordinates of a point $x_0 \in D$ through which, as an initial value, the existence of a solution $\mathbf{x}(t)$ of \eqref{eq13} is verified. Also it is useful in further constructions to provide a sufficient condition for continuity of $\varphi$ or, equivalently, that of $\varphi(t,x_0)$ jointly in $(t,x_0)$:  

\begin{theorem} \thlabel{T24}
If for every $x_0 \in D$ there exists a unique solution $t \rightarrow \varphi(t,x_0)$ of \eqref{eq13} defined on $\mathbb{R}$, then the function $\varphi:\mathbb{R} \times D \rightarrow D$ is continuous.
\end{theorem} 

\begin{proof}
Take $(\tau,x_0) \in \mathbb{R} \times D$ arbitrarily. By hypothesis there is a unique solution $t \rightarrow \varphi(t,x_0)$ on $\mathbb{R}$. So by theorem 3.4 in chapter 1 of \cite{Hale}, $\varphi$ is continuous at $(\tau,x_0)$. Therefore $\varphi$ is continuous on $\mathbb{R} \times D$.
\end{proof}

Assume further that for each $\xi\in D$ a unique solution $\mathbf{x} : \mathbb{R} \rightarrow D$ of \eqref{eq13} exists such that $\mathbf{x}(0)=\xi$. Then the flow 
\begin{equation*}
    \phi: D\times \mathbb{R} \rightarrow D
\end{equation*}
induced by \eqref{eq13} is defined by $\phi(\xi,t)=\mathbf{x}(t)$, where $\mathbf{x}$ is the solution to \eqref{eq13} with initial condition $\mathbf{x}(0)=\xi$. Note that $\phi$ naturally defines a dynamical system $(D,\mathbb{R}, \phi)$, since for every $y \in D$ and $t_1,t_2 \in \mathbb{R}$ we have 
\begin{align*}
   & \phi(y,0)=y,&\text{[identity property]} \\
    &\phi(\phi(y,t_1),t_2)=\phi(y,t_1+t_2),&\text{[group property]}
\end{align*}
and $\phi$ is continuous. The identity property is obtained by the definition of $\phi$; and the group property by the fact that for all $t_1 \in \mathbb{R}$, the functions $t\rightarrow \phi(\phi(y,t_1),t)$ and $t \rightarrow \phi(y,t_1+t)$ satisfy \eqref{eq13}, and are equal at $t=0$, which implies the equality of the two functions by the uniqueness assumption. Also the continuity property of $\phi$ is obtained by \thref{T24}, as $\phi(x_0,t)=\psi(t,x_0)$ for every $(t,x_0) \in \mathbb{R} \times D$. \\

Now, we are at the position to provide the main theorem regarding asymptotic stability of compact subsets of $D$.\\

\begin{theorem}
Assume for the system of differential equations \eqref{eq13} there exists a unique forward solution from each $\xi \in D$. Let $M \subset D$ be a non-empty compact set. If there exists a continuously differentiable real valued function $L$ defined on a neighborhood $N$ of $M$ such that\\

\noindent (i) $L(x)=0$, if $x \in M$, and $L(x)>0$ if $x \notin M$,\\
(ii) $\nabla L(x)^{\top} V(x)<0$, for all $x \notin M$.\\

\noindent Then $M$ is (Lyapunov) asymptotically stable.
\end{theorem}

\begin{proof}
Due to the facts that $L$ is continuously differentiable and $V$ is continuous, the function $\psi: N \rightarrow \mathbb{R}$, with $$\psi(x)=\nabla L(x)^{\top} V(x),$$
is continuous. First of all, note that $\psi(M^{\mathrm{o}})=\{0\}$, since $L$ being continuously differentiable and $L(M)=\{0\}$ imply $\nabla L(M^{\mathrm{o}})=\{0\}$. Furthermore, continuity of $\psi$ will imply $\psi(\partial{M})\subset (-\infty, 0]$. Hence, $\psi\leq0$ by hypothesis (ii). Second, 
\begin{align*}
    \psi(\phi(x,t))&=\nabla L(\phi(x,t))^{\top} V(\phi(x,t)) \notag\\
                                &=\frac{\partial L(\phi(x,t))}{\partial t}.
\end{align*}
So for every $x \notin M$ and $t>0$ we have 
\begin{align}
    L(\phi(x,t))&=L(x)+\int_{0}^{t} \frac{\partial L(\phi(x,t'))}{\partial t'} dt' \notag\\
                &=L(x)+\int_{0}^{t} \psi(\phi(x,t'))dt' \notag\\
                &<L(x), \label{eq14}
\end{align}
where inequality \eqref{eq14} is true since, as we proved, $\psi \leq 0$ in general, and $\psi(\phi(x,0))=\psi(x)<0$ by hypothesis (ii) specifically. Therefore, \thref{T18} immediately concludes $M$ is (Lyapnunov) asymptotically stable.
\end{proof} 

\pagebreak

\end{document}